\documentclass[11pt]{article}
\usepackage{latexsym,amssymb,amsmath,ulem}
\usepackage{color}
\usepackage{graphicx}
\usepackage{amsfonts}

\numberwithin{equation}{section}
\makeatletter

\ifx\pdfoutput\relax\let\pdfoutput=\undefined\fi
\newcount\msipdfoutput
\ifx\pdfoutput\undefined
\else
 \ifcase\pdfoutput
 \else
    \ifx\paperwidth\undefined
    \else
      \ifdim\paperheight=0pt\relax
      \else
        \pdfpageheight\paperheight
      \fi
      \ifdim\paperwidth=0pt\relax
      \else
        \pdfpagewidth\paperwidth
      \fi
    \fi
  \fi
\fi

\newcount\@hour\newcount\@minute\chardef\@x10\chardef\@xv60
\def\tcitime{
\def\@time{%
  \@minute\time\@hour\@minute\divide\@hour\@xv
  \ifnum\@hour<\@x 0\fi\the\@hour:%
  \multiply\@hour\@xv\advance\@minute-\@hour
  \ifnum\@minute<\@x 0\fi\the\@minute
  }}%

\def\x@hyperref#1#2#3{%
   \catcode`\~ = 12
   \catcode`\$ = 12
   \catcode`\_ = 12
   \catcode`\# = 12
   \catcode`\& = 12
   \catcode`\% = 12
   \y@hyperref{#1}{#2}{#3}%
}

\def\y@hyperref#1#2#3#4{%
   #2\ref{#4}#3
   \catcode`\~ = 13
   \catcode`\$ = 3
   \catcode`\_ = 8
   \catcode`\# = 6
   \catcode`\& = 4
   \catcode`\% = 14
}

\@ifundefined{hyperref}{\let\hyperref\x@hyperref}{}
\@ifundefined{msihyperref}{\let\msihyperref\x@hyperref}{}

\@ifundefined{qExtProgCall}{\def\qExtProgCall#1#2#3#4#5#6{\relax}}{}
\def\QCTOpt[#1]#2{%
  \def\QCTOptB{#1}
  \def\QCTOptA{#2}
}
\def\QCTNOpt#1{%
  \def\QCTOptA{#1}
  \let\QCTOptB\empty
}
\def\Qct{%
  \@ifnextchar[{%
    \QCTOpt}{\QCTNOpt}
}
\def\QCBOpt[#1]#2{%
  \def\QCBOptB{#1}%
  \def\QCBOptA{#2}%
}
\def\QCBNOpt#1{%
  \def\QCBOptA{#1}%
  \let\QCBOptB\empty
}
\def\Qcb{%
  \@ifnextchar[{%
    \QCBOpt}{\QCBNOpt}%
}
\def\PrepCapArgs{%
  \ifx\QCBOptA\empty
    \ifx\QCTOptA\empty
      {}%
    \else
      \ifx\QCTOptB\empty
        {\QCTOptA}%
      \else
        [\QCTOptB]{\QCTOptA}%
      \fi
    \fi
  \else
    \ifx\QCBOptA\empty
      {}%
    \else
      \ifx\QCBOptB\empty
        {\QCBOptA}%
      \else
        [\QCBOptB]{\QCBOptA}%
      \fi
    \fi
  \fi
}
\newcount\GRAPHICSTYPE
\GRAPHICSTYPE=\z@
\def\GRAPHICSPS#1{%
 \ifcase\GRAPHICSTYPE%\GRAPHICSTYPE=0
   \special{ps: #1}%
 \or%\GRAPHICSTYPE=1
   \special{language "PS", include "#1"}%
 \fi
}%
\def\graffile#1#2#3#4{%
    \bgroup
	   \@inlabelfalse
       \leavevmode
       \@ifundefined{bbl@deactivate}{\def~{\string~}}{\activesoff}%
        \raise -#4 \BOXTHEFRAME{%
           \hbox to #2{\raise #3\hbox to #2{\null #1\hfil}}}%
    \egroup
}%
\def\draftbox#1#2#3#4{%
 \leavevmode\raise -#4 \hbox{%
  \frame{\rlap{\protect\tiny #1}\hbox to #2%
   {\vrule height#3 width\z@ depth\z@\hfil}%
  }%
 }%
}%
\newcount\@msidraft
\@msidraft=\z@
\let\nographics=\@msidraft
\newif\ifwasdraft
\wasdraftfalse

\def\GRAPHIC#1#2#3#4#5{%
   \ifnum\@msidraft=\@ne\draftbox{#2}{#3}{#4}{#5}%
   \else\graffile{#1}{#3}{#4}{#5}%
   \fi
}
\def\addtoLaTeXparams#1{%
    \edef\LaTeXparams{\LaTeXparams #1}}%
\newif\ifBoxFrame \BoxFramefalse
\newif\ifOverFrame \OverFramefalse
\newif\ifUnderFrame \UnderFramefalse

\def\BOXTHEFRAME#1{%
   \hbox{%
      \ifBoxFrame
         \frame{#1}%
      \else
         {#1}%
      \fi
   }%
}

\def\doFRAMEparams#1{\BoxFramefalse\OverFramefalse\UnderFramefalse\readFRAMEparams#1\end}%
\def\readFRAMEparams#1{%
 \ifx#1\end%
  \let\next=\relax
  \else
  \ifx#1i\dispkind=\z@\fi
  \ifx#1d\dispkind=\@ne\fi
  \ifx#1f\dispkind=\tw@\fi
  \ifx#1t\addtoLaTeXparams{t}\fi
  \ifx#1b\addtoLaTeXparams{b}\fi
  \ifx#1p\addtoLaTeXparams{p}\fi
  \ifx#1h\addtoLaTeXparams{h}\fi
  \ifx#1X\BoxFrametrue\fi
  \ifx#1O\OverFrametrue\fi
  \ifx#1U\UnderFrametrue\fi
  \ifx#1w
    \ifnum\@msidraft=1\wasdrafttrue\else\wasdraftfalse\fi
    \@msidraft=\@ne
  \fi
  \let\next=\readFRAMEparams
  \fi
 \next
 }%
\def\IFRAME#1#2#3#4#5#6{%
      \bgroup
      \let\QCTOptA\empty
      \let\QCTOptB\empty
      \let\QCBOptA\empty
      \let\QCBOptB\empty
      #6%
      \parindent=0pt
      \leftskip=0pt
      \rightskip=0pt
      \setbox0=\hbox{\QCBOptA}%
      \@tempdima=#1\relax
      \ifOverFrame
          \typeout{This is not implemented yet}%
          \show\HELP
      \else
         \ifdim\wd0>\@tempdima
            \advance\@tempdima by \@tempdima
            \ifdim\wd0 >\@tempdima
               \setbox1 =\vbox{%
                  \unskip\hbox to \@tempdima{\hfill\GRAPHIC{#5}{#4}{#1}{#2}{#3}\hfill}%
                  \unskip\hbox to \@tempdima{\parbox[b]{\@tempdima}{\QCBOptA}}%
               }%
               \wd1=\@tempdima
            \else
               \textwidth=\wd0
               \setbox1 =\vbox{%
                 \noindent\hbox to \wd0{\hfill\GRAPHIC{#5}{#4}{#1}{#2}{#3}\hfill}\\%
                 \noindent\hbox{\QCBOptA}%
               }%
               \wd1=\wd0
            \fi
         \else
            \ifdim\wd0>0pt
              \hsize=\@tempdima
              \setbox1=\vbox{%
                \unskip\GRAPHIC{#5}{#4}{#1}{#2}{0pt}%
                \break
                \unskip\hbox to \@tempdima{\hfill \QCBOptA\hfill}%
              }%
              \wd1=\@tempdima
           \else
              \hsize=\@tempdima
              \setbox1=\vbox{%
                \unskip\GRAPHIC{#5}{#4}{#1}{#2}{0pt}%
              }%
              \wd1=\@tempdima
           \fi
         \fi
         \@tempdimb=\ht1
         \advance\@tempdimb by -#2
         \advance\@tempdimb by #3
         \leavevmode
         \raise -\@tempdimb \hbox{\box1}%
      \fi
      \egroup%
}%
\def\DFRAME#1#2#3#4#5{%
  \vspace\topsep
  \hfil\break
  \bgroup
     \leftskip\@flushglue
	 \rightskip\@flushglue
	 \parindent\z@
	 \parfillskip\z@skip
     \let\QCTOptA\empty
     \let\QCTOptB\empty
     \let\QCBOptA\empty
     \let\QCBOptB\empty
	 \vbox\bgroup
        \ifOverFrame
           #5\QCTOptA\par
        \fi
        \GRAPHIC{#4}{#3}{#1}{#2}{\z@}%
        \ifUnderFrame
           \break#5\QCBOptA
        \fi
	 \egroup
  \egroup
  \vspace\topsep
  \break
}%
\def\FFRAME#1#2#3#4#5#6#7{%
  \@ifundefined{floatstyle}
    {%floatstyle undefined (and float.sty not present), no change
     \begin{figure}[#1]%
    }
    {%floatstyle DEFINED
	 \ifx#1h%Only the h parameter, change to H
      \begin{figure}[H]%
	 \else
      \begin{figure}[#1]%
	 \fi
	}
  \let\QCTOptA\empty
  \let\QCTOptB\empty
  \let\QCBOptA\empty
  \let\QCBOptB\empty
  \ifOverFrame
    #4
    \ifx\QCTOptA\empty
    \else
      \ifx\QCTOptB\empty
        \caption{\QCTOptA}%
      \else
        \caption[\QCTOptB]{\QCTOptA}%
      \fi
    \fi
    \ifUnderFrame\else
      \label{#5}%
    \fi
  \else
    \UnderFrametrue%
  \fi
  \begin{center}\GRAPHIC{#7}{#6}{#2}{#3}{\z@}\end{center}%
  \ifUnderFrame
    #4
    \ifx\QCBOptA\empty
      \caption{}%
    \else
      \ifx\QCBOptB\empty
        \caption{\QCBOptA}%
      \else
        \caption[\QCBOptB]{\QCBOptA}%
      \fi
    \fi
    \label{#5}%
  \fi
  \end{figure}%
 }%
\newcount\dispkind%

\def\makeactives{
  \catcode`\"=\active
  \catcode`\;=\active
  \catcode`\:=\active
  \catcode`\'=\active
  \catcode`\~=\active
}
\bgroup
   \makeactives
   \gdef\activesoff{%
      \def"{\string"}%
      \def;{\string;}%
      \def:{\string:}%
      \def'{\string'}%
      \def~{\string~}%
    }
\egroup

\def\FRAME#1#2#3#4#5#6#7#8{%
 \bgroup
 \ifnum\@msidraft=\@ne
   \wasdrafttrue
 \else
   \wasdraftfalse%
 \fi
 \def\LaTeXparams{}%
 \dispkind=\z@
 \def\LaTeXparams{}%
 \doFRAMEparams{#1}%
 \ifnum\dispkind=\z@\IFRAME{#2}{#3}{#4}{#7}{#8}{#5}\else
  \ifnum\dispkind=\@ne\DFRAME{#2}{#3}{#7}{#8}{#5}\else
   \ifnum\dispkind=\tw@
    \edef\@tempa{\noexpand\FFRAME{\LaTeXparams}}%
    \@tempa{#2}{#3}{#5}{#6}{#7}{#8}%
    \fi
   \fi
  \fi
  \ifwasdraft\@msidraft=1\else\@msidraft=0\fi{}%
  \egroup
 }%
\def\TEXUX#1{"texux"}

\def\limfunc#1{\mathop{\rm #1}}%
\def\func#1{\mathop{\rm #1}\nolimits}%
\long\def\QQQ#1#2{%
     \long\expandafter\def\csname#1\endcsname{#2}}%
\@ifundefined{QTP}{\def\QTP#1{}}{}
\@ifundefined{QEXCLUDE}{\def\QEXCLUDE#1{}}{}
\@ifundefined{Qlb}{}{}
\@ifundefined{Qlt}{}{}
\long\def\QQA#1#2{}%
\def\QTR#1#2{{\csname#1\endcsname {#2}}}%
\def\EXPAND#1[#2]#3{}%
\def\NOEXPAND#1[#2]#3{}%
\def\LaTeXparent#1{}%
\def\ChildStyles#1{}%
\def\ChildDefaults#1{}%
\def\QTagDef#1#2#3{}%

\@ifundefined{correctchoice}{}{}
\@ifundefined{HTML}{\def\HTML#1{\relax}}{}
\@ifundefined{TCIIcon}{\def\TCIIcon#1#2#3#4{\relax}}{}
\if@compatibility
\else
  \providecommand{\UNICODE}[2][]{\protect\rule{.1in}{.1in}}
  \providecommand{\U}[1]{\protect\rule{.1in}{.1in}}
  
\fi

\@ifundefined{lambdabar}{
      
   }{}

\@ifundefined{StyleEditBeginDoc}{}{}
\def\QQfnmark#1{\footnotemark}

\@ifundefined{TCIMAKEINDEX}{}{\makeindex}%
\@ifundefined{abstract}{%
 \def\abstract{%
  \if@twocolumn
   \section*{Abstract (Not appropriate in this style!)}%
   \else \small
   \begin{center}{\bf Abstract\vspace{-.5em}\vspace{\z@}}\end{center}%
   \quotation
   \fi
  }%
 }{%
 }%
\@ifundefined{endabstract}{\def\endabstract
  {\if@twocolumn\else\endquotation\fi}}{}%
\@ifundefined{maketitle}{\def\maketitle#1{}}{}%
\@ifundefined{affiliation}{\def\affiliation#1{}}{}%
\@ifundefined{proof}{}{}%
\@ifundefined{endproof}{}{}%
\@ifundefined{newfield}{\def\newfield#1#2{}}{}%
\@ifundefined{chapter}{\def\chapter#1{\par(Chapter head:)#1\par }%
 \newcount\c@chapter}{}%
\@ifundefined{part}{\def\part#1{\par(Part head:)#1\par }}{}%
\@ifundefined{section}{\def\section#1{\par(Section head:)#1\par }}{}%
\@ifundefined{subsection}{\def\subsection#1%
 {\par(Subsection head:)#1\par }}{}%
\@ifundefined{subsubsection}{\def\subsubsection#1%
 {\par(Subsubsection head:)#1\par }}{}%
\@ifundefined{paragraph}{\def\paragraph#1%
 {\par(Subsubsubsection head:)#1\par }}{}%
\@ifundefined{subparagraph}{\def\subparagraph#1%
 {\par(Subsubsubsubsection head:)#1\par }}{}%
\@ifundefined{therefore}{}{}%
\@ifundefined{backepsilon}{}{}%
\@ifundefined{yen}{}{}%
\@ifundefined{registered}{%
   \def\registered{\relax\ifmmode{}\r@gistered
                    \else$\m@th\r@gistered$\fi}%
 \def\r@gistered{^{\ooalign
  {\hfil\raise.07ex\hbox{$\scriptstyle\rm\text{R}$}\hfil\crcr
  \mathhexbox20D}}}}{}%
\@ifundefined{Eth}{}{}%
\@ifundefined{eth}{}{}%
\@ifundefined{Thorn}{}{}%
\@ifundefined{thorn}{}{}%
\@ifundefined{degree}{}{}%
\newdimen\theight
\@ifundefined{Column}{\def\Column{%
 \vadjust{\setbox\z@=\hbox{\scriptsize\quad\quad tcol}%
  \theight=\ht\z@\advance\theight by \dp\z@\advance\theight by \lineskip
  \kern -\theight \vbox to \theight{%
   \rightline{\rlap{\box\z@}}%
   \vss
   }%
  }%
 }}{}%
\@ifundefined{qed}{\def\qed{%
 \ifhmode\unskip\nobreak\fi\ifmmode\ifinner\else\hskip5\p@\fi\fi
 \hbox{\hskip5\p@\vrule width4\p@ height6\p@ depth1.5\p@\hskip\p@}%
 }}{}%
\@ifundefined{cents}{}{}%
\@ifundefined{tciLaplace}{}{}%
\@ifundefined{tciFourier}{}{}%
\@ifundefined{textcurrency}{}{}%
\@ifundefined{texteuro}{}{}%
\@ifundefined{euro}{}{}%
\@ifundefined{textfranc}{}{}%
\@ifundefined{textlira}{}{}%
\@ifundefined{textpeseta}{}{}%
\@ifundefined{miss}{\def\miss{\hbox{\vrule height2\p@ width 2\p@ depth\z@}}}{}%
\@ifundefined{vvert}{}{}%  %always translated to \left| or \right|
\@ifundefined{tcol}{\def\tcol#1{{\baselineskip=6\p@ \vcenter{#1}} \Column}}{}%
\@ifundefined{dB}{}{}%        %dummy entry in column
\@ifundefined{mB}{}{}%   %column entry
\@ifundefined{nB}{}{}%     %column entry (not math)
\@ifundefined{note}{}{}%
\def\newfmtname{LaTeX2e}
\ifx\fmtname\newfmtname
  \DeclareOldFontCommand{\rm}{\normalfont\rmfamily}{\mathrm}
  \DeclareOldFontCommand{\sf}{\normalfont\sffamily}{\mathsf}
  \DeclareOldFontCommand{\tt}{\normalfont\ttfamily}{\mathtt}
  \DeclareOldFontCommand{\bf}{\normalfont\bfseries}{\mathbf}
  \DeclareOldFontCommand{\it}{\normalfont\itshape}{\mathit}
  \DeclareOldFontCommand{\sl}{\normalfont\slshape}{\@nomath\sl}
  \DeclareOldFontCommand{\sc}{\normalfont\scshape}{\@nomath\sc}
\fi

\def\alpha{{\Greekmath 010B}}%
\def\beta{{\Greekmath 010C}}%
\def\gamma{{\Greekmath 010D}}%
\def\delta{{\Greekmath 010E}}%
\def\epsilon{{\Greekmath 010F}}%
\def\zeta{{\Greekmath 0110}}%
\def\eta{{\Greekmath 0111}}%
\def\theta{{\Greekmath 0112}}%
\def\iota{{\Greekmath 0113}}%
\def\kappa{{\Greekmath 0114}}%
\def\lambda{{\Greekmath 0115}}%
\def\mu{{\Greekmath 0116}}%
\def\nu{{\Greekmath 0117}}%
\def\xi{{\Greekmath 0118}}%
\def\pi{{\Greekmath 0119}}%
\def\rho{{\Greekmath 011A}}%
\def\sigma{{\Greekmath 011B}}%
\def\tau{{\Greekmath 011C}}%
\def\upsilon{{\Greekmath 011D}}%
\def\phi{{\Greekmath 011E}}%
\def\chi{{\Greekmath 011F}}%
\def\psi{{\Greekmath 0120}}%
\def\omega{{\Greekmath 0121}}%
\def\varepsilon{{\Greekmath 0122}}%
\def\vartheta{{\Greekmath 0123}}%
\def\varpi{{\Greekmath 0124}}%
\def\varrho{{\Greekmath 0125}}%
\def\varsigma{{\Greekmath 0126}}%
\def\varphi{{\Greekmath 0127}}%

\def\nabla{{\Greekmath 0272}}
\def\FindBoldGroup{%
   {\setbox0=\hbox{$\mathbf{x\global\edef\theboldgroup{\the\mathgroup}}$}}%
}

\def\Greekmath#1#2#3#4{%
    \if@compatibility
        \ifnum\mathgroup=\symbold
           \mathchoice{\mbox{\boldmath$\displaystyle\mathchar"#1#2#3#4$}}%
                      {\mbox{\boldmath$\textstyle\mathchar"#1#2#3#4$}}%
                      {\mbox{\boldmath$\scriptstyle\mathchar"#1#2#3#4$}}%
                      {\mbox{\boldmath$\scriptscriptstyle\mathchar"#1#2#3#4$}}%
        \else
           \mathchar"#1#2#3#4%
        \fi
    \else
        \FindBoldGroup
        \ifnum\mathgroup=\theboldgroup % For 2e
           \mathchoice{\mbox{\boldmath$\displaystyle\mathchar"#1#2#3#4$}}%
                      {\mbox{\boldmath$\textstyle\mathchar"#1#2#3#4$}}%
                      {\mbox{\boldmath$\scriptstyle\mathchar"#1#2#3#4$}}%
                      {\mbox{\boldmath$\scriptscriptstyle\mathchar"#1#2#3#4$}}%
        \else
           \mathchar"#1#2#3#4%
        \fi     	
	  \fi}

\newif\ifGreekBold  \GreekBoldfalse
\let\SAVEPBF=\pbf
\def\pbf{\GreekBoldtrue\SAVEPBF}%

\@ifundefined{mathletters}{%
  \newcounter{equationnumber}
  \def\mathletters{%
     \addtocounter{equation}{1}
     \edef\@currentlabel{\theequation}%
     \setcounter{equationnumber}{\c@equation}
     \setcounter{equation}{0}%
     \edef\theequation{\@currentlabel\noexpand\alph{equation}}%
  }
  
}{}

\@ifundefined{BibTeX}{%
    \def\BibTeX{{\rm B\kern-.05em{\sc i\kern-.025em b}\kern-.08em
                 T\kern-.1667em\lower.7ex\hbox{E}\kern-.125emX}}}{}%
\@ifundefined{AmS}%
    {\def\AmS{{\protect\usefont{OMS}{cmsy}{m}{n}%
                A\kern-.1667em\lower.5ex\hbox{M}\kern-.125emS}}}{}%
\@ifundefined{AmSTeX}{}{}%
\def\@@eqncr{\let\@tempa\relax
    \ifcase\@eqcnt \def\@tempa{& & &}\or \def\@tempa{& &}%
      \else \def\@tempa{&}\fi
     \@tempa
     \if@eqnsw
        \iftag@
           \@taggnum
        \else
           \@eqnnum\stepcounter{equation}%
        \fi
     \fi
     \global\tag@false
     \global\@eqnswtrue
     \global\@eqcnt\z@\cr}

\def\TCItag{\@ifnextchar*{\@TCItagstar}{\@TCItag}}
\def\@TCItag#1{%
    \global\tag@true
    \global\def\@taggnum{(#1)}%
    \global\def\@currentlabel{#1}}
\def\@TCItagstar*#1{%
    \global\tag@true
    \global\def\@taggnum{#1}%
    \global\def\@currentlabel{#1}}
\def\tint{\msi@int\textstyle\int}%
\def\tiint{\msi@int\textstyle\iint}%
\def\tiiint{\msi@int\textstyle\iiint}%
\def\tiiiint{\msi@int\textstyle\iiiint}%
\def\tidotsint{\msi@int\textstyle\idotsint}%
\def\toint{\msi@int\textstyle\oint}%

\newtoks\temptoksa
\newtoks\temptoksb
\newtoks\temptoksc

\def\msi@int#1#2{%
 \def\@temp{{#1#2\the\temptoksc_{\the\temptoksa}^{\the\temptoksb}}}%
 \futurelet\@nextcs
 \@int
}

\def\@int{%
   \ifx\@nextcs\limits
      \typeout{Found limits}%
      \temptoksc={\limits}%
	  \let\@next\@intgobble%
   \else\ifx\@nextcs\nolimits
      \typeout{Found nolimits}%
      \temptoksc={\nolimits}%
	  \let\@next\@intgobble%
   \else
      \typeout{Did not find limits or no limits}%
      \temptoksc={}%
      \let\@next\msi@limits%
   \fi\fi
   \@next
}%

\def\@intgobble#1{%
   \typeout{arg is #1}%
   \msi@limits
}

\def\msi@limits{%
   \temptoksa={}%
   \temptoksb={}%
   \@ifnextchar_{\@limitsa}{\@limitsb}%
}

\def\@limitsa_#1{%
   \temptoksa={#1}%
   \@ifnextchar^{\@limitsc}{\@temp}%
}

\def\@limitsb{%
   \@ifnextchar^{\@limitsc}{\@temp}%
}

\def\@limitsc^#1{%
   \temptoksb={#1}%
   \@ifnextchar_{\@limitsd}{\@temp}%
}

\def\@limitsd_#1{%
   \temptoksa={#1}%
   \@temp
}

\def\dint{\msi@int\displaystyle\int}%
\def\diint{\msi@int\displaystyle\iint}%
\def\diiint{\msi@int\displaystyle\iiint}%
\def\diiiint{\msi@int\displaystyle\iiiint}%
\def\didotsint{\msi@int\displaystyle\idotsint}%
\def\doint{\msi@int\displaystyle\oint}%

\if@compatibility\else
  \RequirePackage{amsmath}
\fi

\let\DOTSI\relax
\def\RIfM@{\relax\ifmmode}%
\def\FN@{\futurelet\next}%
\newcount\intno@
\def\iint{\DOTSI\intno@\tw@\FN@\ints@}%
\def\iiint{\DOTSI\intno@\thr@@\FN@\ints@}%
\def\iiiint{\DOTSI\intno@4 \FN@\ints@}%
\def\idotsint{\DOTSI\intno@\z@\FN@\ints@}%
\def\ints@{\findlimits@\ints@@}%
\newif\iflimtoken@
\newif\iflimits@
\def\findlimits@{\limtoken@true\ifx\next\limits\limits@true
 \else\ifx\next\nolimits\limits@false\else
 \limtoken@false\ifx\ilimits@\nolimits\limits@false\else
 \ifinner\limits@false\else\limits@true\fi\fi\fi\fi}%
\def\multint@{\int\ifnum\intno@=\z@\intdots@                          %1
 \else\intkern@\fi                                                    %2
 \ifnum\intno@>\tw@\int\intkern@\fi                                   %3
 \ifnum\intno@>\thr@@\int\intkern@\fi                                 %4
 \int}%                                                               %5
\def\multintlimits@{\intop\ifnum\intno@=\z@\intdots@\else\intkern@\fi
 \ifnum\intno@>\tw@\intop\intkern@\fi
 \ifnum\intno@>\thr@@\intop\intkern@\fi\intop}%
\def\intic@{%
    \mathchoice{\hskip.5em}{\hskip.4em}{\hskip.4em}{\hskip.4em}}%
\def\negintic@{\mathchoice
 {\hskip-.5em}{\hskip-.4em}{\hskip-.4em}{\hskip-.4em}}%
\def\ints@@{\iflimtoken@                                              %1
 \def\ints@@@{\iflimits@\negintic@
   \mathop{\intic@\multintlimits@}\limits                             %2
  \else\multint@\nolimits\fi                                          %3
  \eat@}%                                                             %4
 \else                                                                %5
 \def\ints@@@{\iflimits@\negintic@
  \mathop{\intic@\multintlimits@}\limits\else
  \multint@\nolimits\fi}\fi\ints@@@}%
\def\intkern@{\mathchoice{\!\!\!}{\!\!}{\!\!}{\!\!}}%
\def\plaincdots@{\mathinner{\cdotp\cdotp\cdotp}}%
\def\intdots@{\mathchoice{\plaincdots@}%
 {{\cdotp}\mkern1.5mu{\cdotp}\mkern1.5mu{\cdotp}}%
 {{\cdotp}\mkern1mu{\cdotp}\mkern1mu{\cdotp}}%
 {{\cdotp}\mkern1mu{\cdotp}\mkern1mu{\cdotp}}}%
\def\RIfM@{\relax\protect\ifmmode}
\def\text{\RIfM@\expandafter\text@\else\expandafter\mbox\fi}
\let\nfss@text\text
\def\text@#1{\mathchoice
   {\textdef@\displaystyle\f@size{#1}}%
   {\textdef@\textstyle\tf@size{\firstchoice@false #1}}%
   {\textdef@\textstyle\sf@size{\firstchoice@false #1}}%
   {\textdef@\textstyle \ssf@size{\firstchoice@false #1}}%
   \glb@settings}

\def\textdef@#1#2#3{\hbox{{%
                    \everymath{#1}%
                    \let\f@size#2\selectfont
                    #3}}}
\newif\iffirstchoice@
\firstchoice@true

\def\tto{\;{\lower 1pt \hbox{$\rightarrow$}}\kern -10pt
\hbox{\raise 2pt \hbox{$\rightarrow$}}\;}

\newtheorem{theorem}{Theorem}[section]
\newtheorem{proposition}{Proposition}[section]
\newtheorem{corollary}{Corollary}[section]
\newtheorem{lemma}{Lemma}[section]
\newtheorem{remark}{Remark}[section]
\newtheorem{example}{Example}[section]
\newtheorem{definition}{Definition}[section]

\numberwithin{equation}{section}

\renewcommand{\theequation}{\thesection.\arabic{equation}}
\normalsize

\@ifundefined{theorem}{\newtheorem{theorem}{Theorem}}{}
\@ifundefined{lemma}{\newtheorem{lemma}[theorem]{Lemma}}{}
\@ifundefined{corollary}{\newtheorem{corollary}[theorem]{Corollary}}{}
\@ifundefined{conjecture}{}{}
\@ifundefined{proposition}{\newtheorem{proposition}[theorem]{Proposition}}{}
\@ifundefined{axiom}{}{}
\@ifundefined{remark}{\newtheorem{remark}{Remark}}{}
\@ifundefined{example}{\newtheorem{example}{Example}}{}
\@ifundefined{exercise}{}{}
\@ifundefined{definition}{}{}

\ifx\pdfoutput\relax\let\pdfoutput=\undefined\fi
\newcount\msipdfoutput
\ifx\pdfoutput\undefined
\else
 \ifcase\pdfoutput
 \else
    \ifx\paperwidth\undefined
    \else
      \ifdim\paperheight=0pt\relax
      \else
        \pdfpageheight\paperheight
      \fi
      \ifdim\paperwidth=0pt\relax
      \else
        \pdfpagewidth\paperwidth
      \fi
    \fi
  \fi
\fi

\makeatother

\begin{document}
\title{\textbf{Around a Farkas type Lemma }}
\date{\ \ \today }
\author{N. Dinh\thanks{%
International University, Vietnam National University - HCMC, and Vietnam
National University - HCM city, Linh Trung ward, Thu Duc district, Ho Chi
Minh city, Vietnam (ndinh@hcmiu.edu.vn) }, \ \ M. A. Goberna\thanks{%
Department of Mathematics, University of Alicante, Alicante, Spain
(mgoberna@ua.es)} \thanks{%
Corresponding author.}, \ M. Volle\thanks{%
Avignon University, LMA EA 2151, Avignon, France
(michel.volle@univ-avignon.fr)} }
\maketitle

\begin{abstract}
The first two authors of this paper asserted in Lemma 4 of "New Farkas-type
constraint qualifications in convex infinite programming" (DOI:
10.1051/cocv:2007027) that a given reverse convex inequality is consequence
of a given convex system satisfying the Farkas-Minkowski constraint
qualification if and only if certain set depending on the data contains a
particular point of the vertical axis. This paper identifies a hidden
assumption in this reverse Farkas lemma which always holds in its
applications to nontrivial optimization problems. Moreover, it shows that
the statement remains valid when the Farkas-Minkowski constraint
qualification fails by replacing the mentioned set by its closure. This
hidden assumption is also characterized in terms of the data. Finally, the
paper provides some applications to convex infinite systems and to convex
infinite optimization problems.
\end{abstract}

%\title{ Around a Farkas type Lemma  }

\textbf{Key words }Convex infinite systems; Farkas Lemma; Convex infinite
optimization

\textbf{Mathematics Subject Classification }Primary 90C25; Secondary 49N15;
46N10

\section{Introduction}

This paper deals with inequality systems of the form
\begin{equation}
\sigma :=\{f_{i}(x)\leq 0,i\in I;\;x\in C\},  \label{1.2}
\end{equation}%
where $I$ is an arbitrary (possibly infinite) index set, $C$ is a nonempty
closed convex subset of a locally convex Hausdorff space $X$ and $\left\{
f_{i}\right\} _{i\in I}\subset \Gamma \left( X\right) ,$ the set of lower
semicontinuous (lsc in brief) convex proper functions on $X.$

The Farkas type methodology in convex optimization is about the following
question regarding a constraint system $\sigma ,$ an objective function $%
f\in \Gamma \left( X\right) $ and a scalar $\alpha \in \mathbb{R}:$

(Q) \textit{When is the relation }$f\left( x\right) \geq \alpha $\textit{\ a
consequence }$\sigma $\textit{? }

The answer to (Q) in \cite{DGLS07} was the so-called reverse Farkas Lemma
\cite[Lemma 4]{DGLS07}, which provided a characterization expressed in terms
of the data ($C,$ $\left\{ f_{i}\right\} _{i\in I},$ $f$ and $\alpha $).
More precisely, under certain constraint qualification on $\sigma $, the
answer consisted in the membership of certain point of the vertical axis
depending on $\alpha $ to the convex cone generated by the epigraphs of
conjugates of the involved functions, $\left\{ f_{i}\right\} _{i\in I}$ and $%
f,$ and the indicator of $C,$ but omitted a hidden assumption (the existence
of a solution of $\sigma $ where $f$ takes a finite value) which fortunately
holds at any nontrivial application to convex optimization. The new version
of the reverse Farkas Lemma in this paper repairs this drawback and extends
the result in two directions: getting rid of constraint qualifications on
the system $\sigma $ and considering arbitrary linear perturbations of the
function $f.$ The readers are referred to the reviews \cite{DJ14} and \cite%
{Jeya08}, and references therein, to put both versions of the reverse Farkas
Lemma in the context of the vast literature on Farkas type lemmas.

The paper is organized as follows. Section 2 introduces the convex analysis
tools to be used along the paper. Section 3 provides a counterexample
showing that \cite[Lemma 4]{DGLS07} may fail when the hidden assumption does
not hold even if the Farkas-Minkowski constraint qualification is satisfied.
Section 4 contains the main results of the paper: Theorem \ref{thm1}, which
makes explicit the hidden assumption and shows that the statement of \cite[%
Lemma 4]{DGLS07} remains valid without constraint qualifications by
replacing the mentioned set by its closure, Theorem \ref{thm1b}, which
consists in a non-asymptotic stable reverse Farkas Lemma, Theorem \ref{thm2}%
, which characterizes the mentioned hidden assumption in terms of the data
for arbitrary convex systems, and Theorem \ref{thm3}, which provides
necessary conditions for the fulfilment of the hidden assumption,\ $\sigma $
being consistent or inconsistent. Finally, Section 5 provides applications
to convex infinite systems and to convex infinite optimization problems.

\section{Preliminaries}

We shall start this section with some necessary definitions and well known
properties. Given a nonempty subset $\QTR{frametitle}{D}$ of a (real
Hausdorff) locally convex topological vector space $X$ with null vector $%
0_{X},$\ we denote by $\limfunc{co}\QTR{frametitle}{D}$ and $\limfunc{cone}%
\QTR{frametitle}{D}$ the convex hull of $\QTR{frametitle}{D}$ and the convex
cone generated by $\QTR{frametitle}{D}\cup \left\{ 0_{X}\right\} ,$
respectively. One has $\limfunc{cone}\QTR{frametitle}{D=}\limfunc{co}\left(
\mathbb{R}_{+}D\right) =\mathbb{R}_{+}\limfunc{co}\left( D\right) .$ The
(negative) polar of a nonempty set $D\subset X$ is $D^{-}:=\left\{ x^{\ast
}\in X^{\ast }:\left\langle x^{\ast },x\right\rangle \leq 0,\forall x\in
D\right\} .$ The closure of a subset $D\subset X$ will be denoted by $%
\limfunc{cl}D.$\ Similarly, the closure of a subset $D$ of the dual space $%
X^{\ast }$ (respectively the product space $X^{\ast }\times \mathbb{R}$)
equipped with the weak$^{\ast }$ topology (respectively the product
topology) will be denoted by $w^{\ast }-\limfunc{cl}D.$ For the sake of
brevity, we also denote both closures of $D$\ by $\overline{D}$ when no
confusion is possible. Given two subsets $D$ and $E$ of $X,$ one has%
\begin{equation}
\overline{\overline{D}+E}=\overline{D+E}.  \label{2.1}
\end{equation}

The recession cone of a nonempty closed convex set $D\subset X$ is given by%
\begin{equation*}
D^{\infty }=\bigcap\limits_{t>0}t\left( D-a\right)
\end{equation*}%
for any $a\in D.$ For any family $\left\{ D_{j}\right\} _{j\in J}$ of closed
convex sets such that $\bigcap\limits_{j\in J}D_{j}\neq \emptyset ,$ one has
\begin{equation*}
\left( \bigcap\limits_{j\in J}D_{j}\right) ^{\infty }=\bigcap\limits_{j\in
J}D_{j}^{\infty }.
\end{equation*}

Given a function $h:X\longrightarrow \overline{\mathbb{R}}:=\mathbb{R\cup
\{\pm \infty \}},$ its domain is $\limfunc{dom}h:=\{x\in X:h(x)<+\infty \},$
its graph is $\func{gph}h:=\{\left( x,h\left( x\right) \right) :x\in
X,h(x)\in \mathbb{R}\},$ its epigraph is $\limfunc{epi}h:=\{\left(
x,r\right) \in X\times \mathbb{R}:h(x)\leq r\}$ and its Fenchel conjugate is
$h^{\ast }:X^{\ast }\longrightarrow \overline{\mathbb{R}}$ such that $%
h^{\ast }(x^{\ast }):=\sup \{\langle x^{\ast },x\rangle -h(x):x\in X\}$ for
any $x^{\ast }\in X^{\ast }$. Given $r\in \mathbb{R}$, we also denote by $%
\left[ h\leq r\right] :=\left\{ x\in X:h(x)\leq r\right\} $\ the lower level
set of $h$ at level $r.$ For a set $D\subset X$, the indicator $\delta _{D}$
of a $D$ is the function that takes value $0$ at $x\in X$\ if $x\in D$ and $%
+\infty $ otherwise. $D$ is a nonempty closed convex subset of $X$ if and
only if $\delta _{D}\in \Gamma \left( X\right) .$

The barrier cone of a nonempty set $D\subset X$ is
\begin{equation*}
\limfunc{barr}D:=\left\{ x^{\ast }\in X^{\ast }:\delta _{D}^{\ast }\left(
x^{\ast }\right) <+\infty \right\} =\func{dom}\delta _{D}^{\ast }.
\end{equation*}

For any pair of nonempty subsets $D$ and $E$ of $X$ one has%
\begin{equation}
\delta _{D+E}^{\ast }=\delta _{D}^{\ast }+\delta _{E}^{\ast }.  \label{2.3}
\end{equation}

The subdifferential of $h$ at $a\in X$ is%
\begin{equation*}
\partial h(a):=\left\{
\begin{array}{ll}
\{x^{\ast }\in X^{\ast }\,:\,h(x)\geq h(a)+\langle x^{\ast },x-a\rangle
,\forall x\in X\}, & \text{if }h(a)\in \mathbb{R}, \\
\emptyset , & \text{else.}%
\end{array}%
\right.
\end{equation*}%
We denote by $N_{D}\left( a\right) :=\partial \delta _{D}(a)$ the normal
cone to $D$ at $a\in D.$

Given $g,h\in \Gamma \left( X\right) $ such that $\limfunc{dom}g\cap
\limfunc{dom}h\neq \emptyset ,$ then
\begin{equation}
\limfunc{epi}\left( g+h\right) ^{\ast }=w^{\ast }-\limfunc{cl}\left(
\limfunc{epi}g^{\ast }+\limfunc{epi}h^{\ast }\right) .  \label{2.2}
\end{equation}%
In particular, when one of the functions $f$ or $g$ is continuous at a point
in the domain of the other, the weak$^{\ast }$ closure (\textquotedblleft $%
w^{\ast }-\limfunc{cl}$") in the right hand side of \eqref{2.2} can be
dropped.

Given $\left\{ h_{i}\right\} _{i\in I}\subset \Gamma \left( X\right) $ with $%
h:=\sup\limits_{i\in I}h_{i}$ proper, one has
\begin{equation}
\limfunc{epi}h^{\ast }=w^{\ast }-\limfunc{cl}\limfunc{co}\left(
\bigcup\limits_{i\in I}\limfunc{epi}h_{i}^{\ast }\right) .  \label{2.4}
\end{equation}

For\ a proper function $g:X\rightarrow \overline{\mathbb{R}}$ one has\textbf{%
\ }
\begin{equation*}
\delta _{\limfunc{epi}g}^{\ast }(x^{\ast },r)=\left\{
\begin{array}{ll}
-rg^{\ast }(-r^{-1}x^{\ast }), & r<0, \\
\delta _{\limfunc{dom}g}^{\ast }(x^{\ast }), & r=0, \\
+\infty , & r>0.%
\end{array}%
\right.
\end{equation*}%
Applying this equation to $g:=h^{\ast }$, for $h\in \Gamma \left( X\right) ,$
one gets the perspective function%
\begin{equation}
\delta _{\limfunc{epi}h^{\ast }}^{\ast }\left( x,r\right) =\left\{
\begin{array}{ll}
-rh\left( -\frac{x}{r}\right) , & \text{if }r<0, \\
\delta _{\func{dom}h^{\ast }}^{\ast }\left( x\right) , & \text{if }r=0, \\
+\infty , & \text{if }r>0.%
\end{array}%
\right.  \label{2.6}
\end{equation}%
{The recession function of $h\in \Gamma (X)$ is $h^{\infty }:X^{\ast
}\longrightarrow \overline{\mathbb{R}}$ such that $\limfunc{epi}(h^{\infty
})=(\limfunc{epi}h)^{\infty }$. } Moreover,
\begin{equation}
h^{\infty }=\delta _{\func{dom}h^{\ast }}^{\ast }=\delta _{\limfunc{epi}%
h^{\ast }}^{\ast }\left( \cdot ,0\right) ,  \label{2.7}
\end{equation}

\begin{equation}
\left[ h^{\infty }\leq 0\right] =\left( \func{dom}h^{\ast }\right) ^{-},
\label{2.8}
\end{equation}%
and
\begin{equation}
\left[ h\leq r\right] \neq \emptyset \Longrightarrow \left[ h^{\infty }\leq 0%
\right] =\left[ h\leq r\right] ^{\infty }.  \label{2.9}
\end{equation}

The next lemma is a well-known consequence of the Hahn-Banach Theorem which
will be used in the sequel.

\begin{lemma}
\label{lem1nw} for each $D\subset X,$ one has%
\begin{equation}
\overline{\limfunc{co}}D:=\limfunc{cl}\limfunc{co}D= \{x\in X:\forall
x^{\ast }\in X^{\ast },\left\langle x^{\ast },x\right\rangle \leq \delta
_{D}^{\ast }\left( x^{\ast }\right) \}.  \label{2.5}
\end{equation}
\end{lemma}

The \textit{solution set } of $\sigma $ is%
\begin{equation*}
A:=\{x\in C:f_{i}(x)\leq 0,\forall i\in I\}=B\cap C,
\end{equation*}%
where
\begin{equation*}
B:=\{x\in X:f_{i}(x)\leq 0,\forall i\in I\}=\bigcap\limits_{i\in I}\left[
f_{i}\leq 0\right] .
\end{equation*}

Denoting by $0_{X^{\ast }}$ the null element of $X^{\ast }$,\ the question
(Q) posed in Section 1 can equivalently be formulated as follows: given $%
\sigma $ as in (\ref{1.2}), $f\in \Gamma \left( X\right) $ and $\alpha \in
\mathbb{R},$ \textit{when}
\begin{equation}
\left( 0_{X^{\ast }},-\alpha \right) \in \limfunc{epi}\left( f+\delta
_{A}\right) ^{\ast }  \label{3.1}
\end{equation}%
\textit{holds? }

The answer to question (Q) {will be discussed in the next sections. It}
involves the so-called \textit{characteristic cone} of $\sigma $ \cite{DGL06}%
:
\begin{equation*}
K:=\limfunc{cone}\left\{ \limfunc{epi}\delta _{C}^{\ast }\cup \bigcup_{%
\QTR{frametitle}{i}\in I}\limfunc{epi}f_{\QTR{frametitle}{i}}^{\ast
}\right\} .
\end{equation*}%
Defining $\Delta :=\bigcup\limits_{\QTR{frametitle}{i}\in \QTR{frametitle}{I}%
}\limfunc{epi}f_{\QTR{frametitle}{i}}^{\ast }$ {and noting that $\limfunc{epi%
}\delta _{C}^{\ast }$ is a convex cone,} the characteristic cone can be
expressed as
\begin{equation}
K=\limfunc{cone}\left\{ \limfunc{epi}\delta _{C}^{\ast }\cup \Delta \right\}
=\limfunc{epi}\delta _{C}^{\ast }+\limfunc{cone}\Delta .  \label{3.2}
\end{equation}%
Lastly, it is worth mentioning that
\begin{equation}
\limfunc{cone}\Delta =\bigcup_{J\in \mathcal{J}}\Big\{\sum_{j\in J}\lambda
_{j}\left( x_{j}^{\ast },r_{j}\right) :\lambda _{j}\geq 0,\left( x_{j}^{\ast
},r_{j}\right) \in \limfunc{epi}f_{j}^{\ast },\forall j\in J\Big\},
\label{2.10}
\end{equation}%
where $\mathcal{J}$ is the collection of all nonempty finite subsets of $I$.

Given $J\in \mathcal{J}$, $\lambda _{j}\geq 0$ for all $j\in J$, we can
consider the finite sum
\begin{equation*}
\Big(\sum_{j\in J}\lambda _{j}f_{j}\Big)(x)=\sum_{j\in J}\lambda
_{j}f_{j}(x),x\in X,
\end{equation*}%
in which the rule $0\times (+\infty )=0$ is adopted.

\section{A counterexample}

Following \cite{DGL06},\ the system $\sigma $ in (\ref{1.2}) is said to be
\textit{Farkas-Minkowski} (FM, in short) if $A\neq \emptyset $ and $K $ is $%
w^{\ast }-$closed. In \cite[Lemma 4]{DGLS07} it is claimed that if $\sigma $
is FM, then%
\begin{equation*}
\text{(\ref{3.1})}\Longleftrightarrow \left( 0_{X^{\ast }},-\alpha \right)
\in w^{\ast }-\limfunc{cl}\left( \limfunc{epi}f^{\ast }+K\right) .
\end{equation*}

Let us give an example showing that under the assumption that $K$ is $%
w^{\ast }-$closed and $A\neq \emptyset ,$ the statement%
\begin{equation*}
\text{(\ref{3.1})}\Longrightarrow \left( 0_{X^{\ast }},-\alpha \right) \in
w^{\ast }-\limfunc{cl}\left( \limfunc{epi}f^{\ast }+K\right)
\end{equation*}%
may fail.

\begin{example}
\label{Example1}Let $X=C=\mathbb{R}^{2},$ $I=\left\{ 1\right\} ,$
\begin{equation*}
f_{1}\left( u,v\right) =\left\{
\begin{array}{ll}
0, & \text{if }\left( u,v\right) \in \mathbb{R}_{+}^{2}\text{ and }v=u+1, \\
+\infty , & \text{else,}%
\end{array}%
\right.
\end{equation*}%
and%
\begin{equation*}
f\left( u,v\right) =\left\{
\begin{array}{ll}
-u, & \text{if }\left( u,v\right) \in \mathbb{R}_{+}^{2}\text{ and }v=u, \\
+\infty , & \text{else.}%
\end{array}%
\right.
\end{equation*}%
We have
\begin{equation*}
A=B=\left[ f_{1}\leq 0\right] =\left\{ \left( 0,1\right) \right\} +\mathbb{R}%
_{+}\left\{ \left( 1,1\right) \right\} ,
\end{equation*}%
$f_{1}=\delta _{A},$ and%
\begin{equation*}
\begin{array}{ll}
K & {=}\limfunc{cone}\left\{ \limfunc{epi}f_{1}^{\ast }\cup \left( \left\{
0_{\mathbb{R}^{2}}\right\} \times \mathbb{R}_{+}\right) \right\} =\limfunc{%
cone}\limfunc{epi}\delta _{A}^{\ast } \\
& =\limfunc{epi}\delta _{A}^{\ast }=\left\{ \left( a,b,r\right) \in \mathbb{R%
}^{3}:a+b\leq 0,b\leq r\right\} ,%
\end{array}%
\end{equation*}%
so that $A\neq \emptyset $ and $\sigma $ is FM.\ Note that $A\cap \func{dom}%
f=\emptyset $ and, consequently,%
\begin{equation*}
\left( f+\delta _{A}\right) ^{\ast }{(0_{\mathbb{R}^{2}})}\leq \alpha
,\forall \alpha \in \mathbb{R},
\end{equation*}%
so that (\ref{3.1}) holds for any $\alpha \in \mathbb{R}.$ Let us prove
that, however,
\begin{equation}
\left( 0_{\mathbb{R}^{2}},-\alpha \right) \notin w^{\ast }-\limfunc{cl}%
\left( \limfunc{epi}f^{\ast }+K\right) ,\forall \alpha \in \mathbb{R}.
\label{eq1nw}
\end{equation}%
{Take $\alpha \in \mathbb{R}$ and note that $\limfunc{epi}f^{\ast }+K$ is
convex. It follows from Lemma \ref{lem1nw} (see (\ref{2.5})) that (\ref%
{eq1nw}) is equivalent to the existence of a triple $\left( a,b,r\right) \in
\mathbb{R}^{3}$ such that
\begin{equation}
-\alpha r>\delta _{\limfunc{epi}f^{\ast }+K}^{\ast }\left( a,b,r\right)
=\delta _{\limfunc{epi}f^{\ast }}^{\ast }\left( a,b,r\right) +\delta
_{K}^{\ast }\left( a,b,r\right)  \label{eq2nw}
\end{equation}%
(the last equation comes from (\ref{2.3})). }

This is the case for $\left( a,b,r\right) =\left( 1,1,0\right).$ One has in
fact $f^{\ast }=\delta _{\left\{ \left( u^{\prime },v^{\prime }\right) \in
\mathbb{R}^{2}:u^{\prime }+v^{\prime }\leq -1\right\} },$%
\begin{equation*}
\delta _{A}^{\ast }\left( u^{\prime },v^{\prime }\right) =\left\{
\begin{array}{ll}
v^{\prime }, & \text{if }u^{\prime },v^{\prime }\leq 0, \\
+\infty , & \text{else,}%
\end{array}%
\right.
\end{equation*}%
\begin{equation*}
\delta _{\limfunc{epi}f^{\ast }}^{\ast }\left( 1,1,0\right)
=\sup\limits_{\left( u^{\prime },v^{\prime }\right) \in \func{dom}f^{\ast
}}\left( u^{\prime }+v^{\prime }\right) =\sup\limits_{u^{\prime }+v^{\prime
}\leq -1}\left( u^{\prime }+v^{\prime }\right) =-1,
\end{equation*}%
and%
\begin{equation*}
\delta _{K}^{\ast }\left( 1,1,0\right) =\sup\limits_{\left( u^{\prime
},v^{\prime }\right) \in \func{dom}\delta _{A}^{\ast }}\left( u^{\prime
}+v^{\prime }\right) =\sup\limits_{u^{\prime }+v^{\prime }\leq 0}\left(
u^{\prime }+v^{\prime }\right) =0,
\end{equation*}
{meaning that \eqref{eq2nw} holds with $\left( a,b,r\right) =\left(
1,1,0\right) .$ }
\end{example}

\section{The reverse Farkas Lemma revisited}

In order to rectify \cite[Lemma 4]{DGLS07} we need three more or less
classical properties.

\begin{lemma}
\label{lemma1}Assume $B\neq \emptyset .$ Then, the following statements hold:%
\newline
$\mathrm{(i)}$ $\overline{\limfunc{cone}}\Delta :=w^{\ast }-\limfunc{cl}%
\limfunc{cone}\Delta =\limfunc{epi}\delta _{B}^{\ast }.$\newline
$\mathrm{(ii)}$ $\overline{\limfunc{barr}}B:=w^{\ast }-\limfunc{cl}\limfunc{%
barr}B=\overline{\limfunc{cone}}\left( \bigcup\limits_{\QTR{frametitle}{i}%
\in I}\func{dom}f_{\QTR{frametitle}{i}}^{\ast }\right) .$
\end{lemma}

\textbf{Proof.} For each $\QTR{frametitle}{i}\in I$ we have $f_{i}\leq
\delta _{B},$ $f_{\QTR{frametitle}{i}}^{\ast }\geq \delta _{B}^{\ast }$ and $%
\limfunc{epi}f_{\QTR{frametitle}{i}}^{\ast }\subset \limfunc{epi}\delta
_{B}^{\ast }.$ Therefore, $\Delta \subset \limfunc{epi}\delta _{B}^{\ast }$
and $\bigcup\limits_{\QTR{frametitle}{i}\in I}\func{dom}f_{\QTR{frametitle}{i%
}}^{\ast }\subset \limfunc{barr}B.$ Since $\limfunc{epi}\delta _{B}^{\ast }$
(respectively $\limfunc{barr}B$) is a $w^{\ast }-$closed convex cone
(respectively a convex cone containing $0_{X^{\ast }}$), it follows that $%
\overline{\limfunc{cone}}\Delta \subset \limfunc{epi}\delta _{B}^{\ast }$
(respectively $\limfunc{cone}\left( \bigcup\limits_{\QTR{frametitle}{i}\in I}%
\func{dom}f_{\QTR{frametitle}{i}}^{\ast }\right) \subset \limfunc{barr}B$).

$\mathrm{(i)}$ Note that $\sup\limits_{\QTR{frametitle}{i}\in
I,t>0}tf_{i}=\delta _{B}.$ Since $\left\{ tf_{i}\right\} _{i\in
I,t>0}\subset \Gamma \left( X\right) $ and $B\neq \emptyset ,$ by (\ref{2.4}%
), we have%
\begin{equation*}
\begin{array}{ll}
\limfunc{epi}\delta _{B}^{\ast } & =\overline{\limfunc{co}}\Big(%
\bigcup\limits_{\QTR{frametitle}{i}\in I,t>0}\limfunc{epi}(tf_{%
\QTR{frametitle}{i}})^{\ast }\Big)=\overline{\limfunc{co}}\Big(%
\bigcup\limits_{\QTR{frametitle}{i}\in I,t>0}t\limfunc{epi}f_{%
\QTR{frametitle}{i}}^{\ast }\Big) \\
& =\overline{\limfunc{co}}\left( \mathbb{R}_{++}\Delta \right) \subset
\overline{\limfunc{co}}\left( \mathbb{R}_{+}\Delta \right) =\overline{%
\limfunc{cone}}\Delta%
\end{array}%
\end{equation*}%
and, finally, $\overline{\limfunc{cone}}\Delta =\limfunc{epi}\delta
_{B}^{\ast }.$

$\mathrm{(ii)}$ Let $L$ be projection of $X^{\ast }\times \mathbb{R}$ onto $%
X^{\ast },$ which is a linear $w^{\ast }-$continuous mapping. One has%
\begin{equation*}
\begin{array}{ll}
\limfunc{barr}B & =L\left( \limfunc{epi}\delta _{B}^{\ast }\right) =L\left(
\overline{\limfunc{cone}}\Delta \right) \subset \overline{L\left( \limfunc{%
cone}\Delta \right) } \\
& =\overline{\limfunc{cone}}L\left( \Delta \right) =\overline{\limfunc{cone}}%
\left( \bigcup\limits_{\QTR{frametitle}{i}\in I}\func{dom}f_{%
\QTR{frametitle}{i}}^{\ast }\right) \subset \overline{\limfunc{barr}}B%
\end{array}%
\end{equation*}%
and, finally, $\overline{\limfunc{barr}}B=\overline{\limfunc{cone}}\left(
\bigcup\limits_{\QTR{frametitle}{i}\in I}\func{dom}f_{\QTR{frametitle}{i}%
}^{\ast }\right) .\hfill $\noindent $\square $

\begin{lemma}
\label{lemma2}Assume that $A\neq \emptyset .$ Then we have
\begin{equation*}
w^{\ast }-\limfunc{cl}K=\limfunc{epi}\delta _{A}^{\ast }.
\end{equation*}
\end{lemma}

\textbf{Proof.} Since $A=B\cap C\neq \emptyset ,$ one has, by (\ref{2.2}),%
\begin{equation*}
\limfunc{epi}\delta _{A}^{\ast }=\limfunc{epi}\left( \delta _{B}+\delta
_{C}\right) ^{\ast }=\overline{\limfunc{epi}\delta _{B}^{\ast }+\limfunc{epi}%
\delta _{C}^{\ast }}
\end{equation*}%
and, by Lemma \ref{lemma1}\textrm{(i)}$,$ $\limfunc{epi}\delta _{B}^{\ast }=%
\overline{\limfunc{cone}}\Delta .$ Consequently, by (\ref{2.1}),%
\begin{equation*}
\limfunc{epi}\delta _{A}^{\ast }=\overline{\overline{\limfunc{cone}}\Delta +%
\limfunc{epi}\delta _{C}^{\ast }}=\overline{\limfunc{cone}\Delta +\limfunc{%
epi}\delta _{C}^{\ast }}=\overline{K}.\ \ \ \ \ \ \ \ \ \ \ \ \ \ \ \qquad
\qquad \qquad \qquad \square
\end{equation*}

Variants of Lemma \ref{lemma2} can be found in \cite[Lemma 3.1]{Jeya03} and
\cite[Lemma 2.1]{JLD06}, where $X$ is $\mathbb{R}^{n}$\ and a Banach space,
respectively.

\begin{lemma}
\label{lemma3}Assume that $A\cap \func{dom}f\neq \emptyset .$ Then we have%
\begin{equation*}
w^{\ast }-\limfunc{cl}\left( \limfunc{epi}f^{\ast }+K\right) =\limfunc{epi}%
\left( f+\delta _{A}\right) ^{\ast }.
\end{equation*}
\end{lemma}

\textbf{Proof.} By (\ref{2.2}), Lemma \ref{lemma2} and (\ref{2.1}),%
\begin{equation*}
\begin{array}{ll}
\limfunc{epi}\left( f+\delta _{A}\right) ^{\ast } & =\overline{\limfunc{epi}%
f^{\ast }+\limfunc{epi}\delta _{A}^{\ast }} \\
& =\overline{\limfunc{epi}f^{\ast }+\overline{K}} \ \ \ \ \  \\
& =\overline{\limfunc{epi}f^{\ast }+K}. \ \ \ \ \ \ \ \ \ \ \ \ \ \ \ \qquad
\qquad\qquad \qquad \qquad\qquad \square%
\end{array}%
\end{equation*}

The next result is an immediate consequence of Lemma \ref{lemma3}. It
rectifies and extends \cite[Lemma 4]{DGLS07}. It characterizes when a
linearly perturbed reverse convex inequality is consequence of $\sigma .$

\begin{theorem}[Characterization of the stable reverse Farkas Lemma]
\label{thm1} Assume that $A\cap \func{dom}f\neq \emptyset .$ For each $%
\left( x^{\ast },s\right) \in X^{\ast }\times \mathbb{R}$ we have%
\begin{equation*}
\left[ \left\{
\begin{array}{c}
f_{i}(x)\leq 0,i\in I \\
\;x\in C%
\end{array}%
\right\} \Longrightarrow f\left( x\right) -\left\langle x^{\ast
},x\right\rangle \geq s\right] \Longleftrightarrow \left( x^{\ast
},-s\right) \in \overline{\limfunc{epi}f^{\ast }+K}.
\end{equation*}
\end{theorem}

\textbf{Proof.} \ Note that the statement between square brackets means that
$\left( x^{\ast },-s\right) \in \limfunc{epi}\left( f+\delta _{A}\right)
^{\ast }$ and the conclusion comes from Lemma \ref{lemma3}.$\hfill $%
\noindent $\hfill $\noindent $\square $

\begin{remark}
\label{rem0} Lemma \ref{lemma3} formula is the key point for the asymptotic
version of stable reverse Farkas lemma described by Theorem \ref{thm1}. The
non-asymptotic stable reverse Farkas lemma is naturally given by the next
Theorem \ref{thm1b}, assuming that $\limfunc{epi}f^{\ast }+K$ is weak$^{\ast
}$-closed. Lemma \ref{lemma4} below gives a sufficient condition for that.
\end{remark}

\begin{lemma}
\label{lemma4} Assume $\sigma $ is FM and $f$ is finite continuous at some
point of $A.$ Then $\limfunc{epi}f^{\ast }+K=\limfunc{epi}\left( f+\delta
_{A}\right) ^{\ast }$ and $\limfunc{epi}f^{\ast }+K$ is $w^{\ast }-$closed.
\end{lemma}

\textbf{Proof.} We have $\limfunc{epi}\left( f+\delta _{A}\right) ^{\ast }=%
\limfunc{epi}f^{\ast }+\limfunc{epi}\delta _{A}^{\ast }$ and, by Lemma \ref%
{lemma2}, $\limfunc{epi}\delta _{A}^{\ast }=\overline{K}.$ Since $\sigma $
is FM we then have%
\begin{equation*}
\limfunc{epi}f^{\ast }+K=\limfunc{epi}f^{\ast }+\overline{K}=\limfunc{epi}%
\left( f+\delta _{A}\right) ^{\ast },
\end{equation*}%
which is $w^{\ast }-$closed.$\hfill $\noindent $\hfill $\noindent $\square $%
\medskip

As a direct consequence of Theorem \ref{thm1} (with $x^{\ast }=0_{X^{\ast }}$%
) and Lemma \ref{lemma4} we get:

\begin{corollary}
\label{corol1}\textrm{(\cite[Lemma 4]{DGLS07}} with the hidden assumption%
\textrm{)} Assume that $A\cap \func{dom}f\neq \emptyset .$ Then $f\left(
x\right) \geq s $ is a consequence of $\sigma $ if and only if $\left(
0_{X^{\ast }},-s\right) \in \overline{\limfunc{epi}f^{\ast }+K},$ where the
latter set can be replaced by $\limfunc{epi}f^{\ast }+K$ whenever $\sigma $
is FM and $f $ is finite continuous at some point of $A.$
\end{corollary}

\begin{theorem}[Non-asymptotic stable reverse Farkas Lemma]
\label{thm1b} Assume that $A\cap \func{dom}f\neq \emptyset $ and $\limfunc{%
epi}f^{\ast }+K$ is weak$^{\ast }$-closed. Then, for each $\left( x^{\ast
},s\right) \in X^{\ast }\times \mathbb{R},$ the following statements are
equivalent:\newline
$\mathrm{(i)}$ $f_{i}(x)\leq 0,i\in I,x\in C\Longrightarrow f\left( x\right)
-\left\langle x^{\ast },x\right\rangle \geq s$. \newline
$\mathrm{(ii)}$ There exist $J\in \mathcal{J}$, $u_{j}^{\ast }\in \limfunc{%
dom}f_{j}^{\ast },$ $\lambda _{j}\geq 0$, for all $j\in J,$ $u^{\ast }\in
\limfunc{dom}f^{\ast },$ and $v^{\ast }\in \limfunc{barr}C$ such that
\begin{equation*}
x^{\ast }=u^{\ast }+v^{\ast }+\sum_{j\in J}\lambda _{j}u_{j}^{\ast }
\end{equation*}%
and%
\begin{equation}
f^{\ast }\left( u^{\ast }\right) +\delta _{C}^{\ast }\left( v^{\ast }\right)
+\sum_{j\in J}\lambda _{j}f_{j}^{\ast }\left( u_{j}^{\ast }\right) \leq -s.
\label{eqnii}
\end{equation}%
\noindent $\mathrm{(iii)}$ There exist $J\in \mathcal{J}$ and $\lambda
_{i}\geq 0$, for all $i\in J,$ such that
\begin{equation*}
f(x)-\langle x^{\ast },x\rangle +\sum_{i\in J}\lambda _{i}f_{i}(x)\geq s,\ \
\forall x\in C.
\end{equation*}
\end{theorem}

\textbf{Proof.} $[\mathrm{(i)}\Longrightarrow \mathrm{(ii)}]$ By Theorem \ref%
{thm1}, \textrm{(i)} is equivalent to $\left( x^{\ast },-s\right) \in
\overline{\limfunc{epi}f^{\ast }+K}.$ Since $\limfunc{epi}f^{\ast }+K$ is
weak$^{\ast }$-closed, one has by (\ref{3.2})
\begin{equation}
\left( x^{\ast },-s\right) \in \limfunc{epi}f^{\ast }+\limfunc{epi}\delta
_{C}^{\ast }+\limfunc{cone}\bigcup_{j\in J}\limfunc{epi}f_{j}^{\ast }.
\label{4.8}
\end{equation}%
It follows from (\ref{4.8}) and (\ref{2.10}) that there exist $J\in \mathcal{%
J}$, $(u_{j}^{\ast },r_{j})\in \limfunc{epi}f_{j}^{\ast }$, for all $j \in J$%
, $(u^{\ast },r)\in \limfunc{epi}f$ and $(v^{\ast },p)\in \limfunc{epi}%
\delta _{C}^{\ast }$ such that%
\begin{equation*}
\left( x^{\ast },-s\right) =(u^{\ast },r)+(v^{\ast },p)+\sum_{j\in J}\lambda
_{j}(u_{j}^{\ast },r_{j}).
\end{equation*}%
We then have $-s\geq f^{\ast }\left( u^{\ast }\right) +\delta _{C}^{\ast } {%
(v^\ast)}+\sum_{j\in J}\lambda _{j}f_{j}^{\ast }( u_{j}^{\ast }) $, with $%
x^{\ast }=u^{\ast }+v^{\ast }+\sum_{j\in J}\lambda _{j}u_{j}^{\ast }$, which
is $\text{(ii)}.$

$[\mathrm{(ii)}\Longrightarrow \mathrm{(iii)}]$ It follows from \eqref{eqnii}
that
\begin{equation*}
f(x)-\langle u^{\ast },x\rangle -\langle v^{\ast },x\rangle +\sum_{j\in
J}\lambda _{j}\Big(f_{j}(x)-\langle u_{j}^{\ast },x\rangle \Big)\geq
s,\forall x\in C,
\end{equation*}%
and, as $x^{\ast }=u^{\ast }+v^{\ast }+\sum_{j\in J}\lambda _{j}u_{j}^{\ast
},$ we get
\begin{equation*}
f(x)-\langle x^{\ast },x\rangle +\sum_{j\in J}\lambda _{j}f_{j}(x)\geq s,\ \
\forall x\in C
\end{equation*}%
and $\mathrm{(iii)}$ follows.

$[\mathrm{(iii)}\Longrightarrow \mathrm{(i)}]$ Note that if $x\in A$ then $%
\sum_{i\in J}\lambda _{i}f_{i}(x)\leq 0$ and it comes from (iii) that $%
f(x)-\langle x^{\ast },x\rangle \geq s,$ which is (i). $\hfill $\noindent $%
\square $

\begin{corollary}
\label{corol0} Assume that $A\cap \func{dom}f\neq \emptyset $ and $\limfunc{%
epi}f^{\ast }+K$ is $w^{\ast }$-closed, and let $s\in \mathbb{R}.$ The
following statements are equivalent:\newline
$\mathrm{(i)}$ $f(x)\geq s$ is a consequence of $\sigma .$ \newline
$\mathrm{(ii)}$ $\left( 0_{X^{\ast }},-s\right) \in \limfunc{epi}f^{\ast
}+K. $ \newline
$\mathrm{(iii)}$ There exist $J\in \mathcal{J}$, $u_{j}^{\ast }\in \limfunc{%
dom}f_{j}^{\ast },$ $\lambda _{j}\geq 0$, for all $j\in J,$ $u^{\ast }\in
\limfunc{dom}f^{\ast },$ and $v^{\ast }\in \limfunc{barr}C$ such that%
\begin{equation*}
u^{\ast }+v^{\ast }+\sum_{j\in J}\lambda _{j}u_{j}^{\ast }=0_{X^{\ast }}
\end{equation*}%
and%
\begin{equation*}
f^{\ast }\left( u^{\ast }\right) +\delta _{C}^{\ast }\left( v^{\ast }\right)
+\sum_{j\in J}\lambda _{j}f_{j}^{\ast }\left( u_{j}^{\ast }\right) \leq -s.
\end{equation*}%
\noindent $\mathrm{(iv)}$ There exist $J\in \mathcal{J}$ and $\lambda
_{j}\geq 0$, for all $j\in J$, such that
\begin{equation*}
f(x)+\sum_{j\in J}\lambda _{j}f_{j}(x)\geq s,\ \ \forall x\in C.
\end{equation*}
\end{corollary}

\textbf{Proof.} It is straightforward consequence of Theorem \ref{thm1} when
$x^{\ast }=0_{X^{\ast }}.\hfill $\noindent $\square $

\begin{remark}
\label{rem3} In \cite[Theorem 2]{DGLS07} is established under the
assumptions that \textquotedblleft $\sigma $ is FM and $\limfunc{epi}f^{\ast
}+\overline{K}$ is $w^{\ast }-$closed\textquotedblright\ and (the hidden
assumption) $A\cap \func{dom}f\neq \emptyset $, that $f\left( x\right) \geq
s $ is a consequence of $\sigma $ if and only if there exist $J\in \mathcal{J%
}$ and $\lambda _{j}\geq 0$, for all $j\in J,$ such that%
\begin{equation*}
f\left( x\right) +\sum_{j\in J}\lambda _{j}f_{j}\left( x\right) \geq
s,\forall x\in C.
\end{equation*}%
\newline
Let us note that \textquotedblleft $\sigma $ is FM and $\limfunc{epi}f^{\ast
}+\overline{K}$ is $w^{\ast }-$closed\textquotedblright\ obviously entails
that $\limfunc{epi}f^{\ast }+K$ is $w^{\ast }-$closed, which is our
assumption in Theorem \ref{thm1b} (without assuming that $\sigma $ is FM).
\end{remark}

Let us now reexamine and characterize the case when $A\cap \func{dom}%
f=\emptyset .$ Note that this amounts to the inconsistency of the system%
\begin{equation*}
\left\{ f_{i}(x)\leq 0,i\in I;x\in C\cap \func{dom}f\right\} ,
\end{equation*}%
in which case the convex set $\func{dom}f$ may be nonclosed. We need some
additional preliminaries. Recall that $L:X^{\ast }\times \mathbb{%
R\longrightarrow }X^{\ast }$ denotes the projection of $X^{\ast }\times
\mathbb{R}$ onto $X^{\ast }$ and define the set%
\begin{equation*}
\Omega :=L\left( \Delta \right) =L\left( \bigcup_{\QTR{frametitle}{i}\in
\QTR{frametitle}{I}}\limfunc{epi}f_{\QTR{frametitle}{i}}^{\ast }\right)
=\bigcup_{\QTR{frametitle}{i}\in \QTR{frametitle}{I}}\func{dom}f_{%
\QTR{frametitle}{i}}^{\ast }.
\end{equation*}%
Since $L$ is linear, reasoning as in the proof of Lemma \ref{lemma1}, we
have $L\left( \limfunc{cone}\Delta \right) =\limfunc{cone}\Omega ,$
\begin{equation*}
\begin{array}{ll}
L\left( \limfunc{epi}f^{\ast }+K\right) & =L\left( \limfunc{epi}f^{\ast }+%
\limfunc{cone}\Delta +\limfunc{epi}\delta _{C}^{\ast }\right) \\
& =\func{dom}f^{\ast }+\limfunc{cone}\Omega +\limfunc{barr}C%
\end{array}%
\end{equation*}%
and, consequently,
\begin{equation*}
\limfunc{epi}f^{\ast }+K\subset \left( \func{dom}f^{\ast }+\limfunc{cone}%
\Omega +\limfunc{barr}C\right) \times \mathbb{R},
\end{equation*}%
as well as
\begin{equation}
\overline{\limfunc{epi}f^{\ast }+K}\subset \overline{\func{dom}f^{\ast }+%
\limfunc{cone}\Omega +\limfunc{barr}C}\times \mathbb{R}.  \label{4.3}
\end{equation}

The dual cone of $\Delta $ is
\begin{equation}
\begin{array}{ll}
\Delta ^{-} & =\left\{ \left( x,r\right) \in X\times \mathbb{R}:\delta
_{\Delta }^{\ast }\left( x,r\right) \leq 0\right\} \\
& =\left( \bigcup\limits_{\QTR{frametitle}{i}\in \QTR{frametitle}{I}}%
\limfunc{epi}f_{\QTR{frametitle}{i}}^{\ast }\right) ^{-} \\
& =\bigcap\limits_{\QTR{frametitle}{i}\in \QTR{frametitle}{I}}\left(
\limfunc{epi}f_{\QTR{frametitle}{i}}^{\ast }\right) ^{-}.%
\end{array}
\label{4.1}
\end{equation}%
Moreover, since $\Delta ^{-}=\left( \limfunc{cone}\Delta \right) ^{-},$
\begin{equation}
\delta _{\limfunc{cone}\Delta }^{\ast }=\delta _{\Delta ^{-}.}  \label{4.2}
\end{equation}

\begin{theorem}[Characterization of the hidden assumption]
\label{thm2}The following statements are equivalent:\newline
$\mathrm{(i)}$ $A\cap \func{dom}f=\emptyset .$\newline
\textrm{(ii)} $\overline{\limfunc{epi}f^{\ast }+K}=\overline{\func{dom}%
f^{\ast }+\limfunc{cone}\Omega +\limfunc{barr}C}\times \mathbb{R}.$
\end{theorem}

\textbf{Proof.} \ We shall prove the equivalence of the negations of \textrm{%
(i)} and \textrm{(i)}$,$ say $\rceil $\textrm{(i)} and $\rceil $\textrm{(i)}$%
,$ respectively.

$\Big[\rceil $\textrm{(i)}$\Longrightarrow \rceil $\textrm{(ii)}$\Big]$ By
Lemma \ref{lemma3} we have $\overline{\limfunc{epi}f^{\ast }+K}=\limfunc{epi}%
\left( f+\delta _{A}\right) ^{\ast },$ with $\left( f+\delta _{A}\right)
^{\ast }$ proper. Therefore $\overline{\limfunc{epi}f^{\ast }+K}$ does not
contain any vertical line. Picking $a^{\ast }\in \func{dom}f^{\ast },$ the
vertical line $\left\{ a^{\ast }\right\} \times \mathbb{R}$ is contained in $%
\overline{\func{dom}f^{\ast }+\limfunc{cone}\Omega +\limfunc{barr}C}\times
\mathbb{R}.$ That means that \textrm{(i)} fails.

$\Big[\rceil \mathrm{(ii)}\Longrightarrow \rceil \mathrm{(i)}\Big]$ By (\ref%
{4.3}), $\rceil $\textrm{(ii)} means that there exists
\begin{equation}
\left( a^{\ast },q\right) \in \overline{\func{dom}f^{\ast }+\limfunc{cone}%
\Omega +\limfunc{barr}C}\times \mathbb{R}  \label{4.4}
\end{equation}%
such that $\left( a^{\ast },q\right) \notin \overline{\limfunc{epi}f^{\ast
}+K}.$ {By Lemma \ref{lem1nw} }there exists $\left( a,p\right) \in X\times
\mathbb{R},$ $\left( a,p\right) \neq \left( 0_{X},0\right) ,$\ such that $%
\left\langle a^{\ast },a\right\rangle +pq>\delta _{\limfunc{epi}f^{\ast
}+K}^{\ast }\left( a,p\right) .$ So, by (\ref{2.3}) and (\ref{3.2}),
\begin{equation}
\left\langle a^{\ast },a\right\rangle +pq>\delta _{\limfunc{epi}f^{\ast
}}^{\ast }\left( a,p\right) +\delta _{\limfunc{cone}\Delta }^{\ast }\left(
a,p\right) +\delta _{\limfunc{epi}\delta _{C}^{\ast }}^{\ast }\left(
a,p\right) .  \label{4.5}
\end{equation}

Note that if $p>0,$ the right-hand side of (\ref{4.5}) is $+\infty $ (by (%
\ref{2.6}) applied to $\delta _{C}$), which is impossible. We thus have $%
p\leq 0.$

Assume now that $p=0.$ Then $a\neq 0_{X};$ by (\ref{2.6}), (\ref{4.5}) and (%
\ref{2.3}), we obtain%
\begin{equation}
\begin{array}{ll}
\left\langle a^{\ast },a\right\rangle & >\sup\limits_{x^{\ast }\in \func{dom}%
f^{\ast }}\left\langle x^{\ast },a\right\rangle +\sup\limits_{\left( x^{\ast
},s\right) \in \limfunc{cone}\Delta }\left\langle x^{\ast },a\right\rangle
+\sup\limits_{x^{\ast }\in \limfunc{barr}C}\left\langle x^{\ast
},a\right\rangle \\
& =\sup\limits_{x^{\ast }\in \func{dom}f^{\ast }}\left\langle x^{\ast
},a\right\rangle +\sup\limits_{x^{\ast }\in \limfunc{cone}\Omega
}\left\langle x^{\ast },a\right\rangle +\sup\limits_{x^{\ast }\in \limfunc{%
barr}C}\left\langle x^{\ast },a\right\rangle \\
& =\delta _{\func{dom}f^{\ast }+\limfunc{cone}\Omega +\limfunc{barr}C}^{\ast
}\left( a\right) .%
\end{array}
\label{4.6}
\end{equation}
Again, by Lemma \ref{lem1nw}, \eqref{4.6} means that
\begin{equation*}
a^{\ast }\notin \overline{\func{dom}f^{\ast }+\limfunc{cone}\Omega +\limfunc{%
barr}C},
\end{equation*}%
in contradiction with (\ref{4.4}).

So, we must have $p<0.$ By (\ref{2.6}), (\ref{4.1}), (\ref{4.2}) and (\ref%
{4.5}), one has%
\begin{equation*}
\left\langle a^{\ast },a\right\rangle +pq>-pf\left( -\frac{a}{p}\right)
+\delta _{\Delta ^{-}}\left( a,p\right) -p\delta _{C}\left( -\frac{a}{p}%
\right) ,
\end{equation*}%
which implies $-\frac{a}{p}\in C\cap \func{dom}f$ and $\left( a,p\right) \in
\Delta ^{-}=\bigcap\limits_{\QTR{frametitle}{i}\in \QTR{frametitle}{I}%
}\left( \limfunc{epi}f_{\QTR{frametitle}{i}}^{\ast }\right) ^{-}.$

Applying again (\ref{2.6}) to $\left( a,p\right) ,$\ for each $i\in I$ we
have%
\begin{equation*}
-pf_{i}\left( -\frac{a}{p}\right) =\delta _{\limfunc{epi}f_{i}^{\ast
}}^{\ast }\left( a,p\right) \leq 0
\end{equation*}%
and, finally,
\begin{equation*}
-\frac{a}{p}\in B\cap C\cap \func{dom}f=A\cap \func{dom}f,
\end{equation*}%
which means that $\rceil $\textrm{(i)} holds.$\hfill $\noindent $\square $

\begin{remark}
\label{rem1}Since $\limfunc{cone}\Omega +\limfunc{barr}C=\limfunc{cone}%
\left( \Omega \cup \limfunc{barr}C\right) $ and $\overline{\limfunc{cone}%
\Omega +\limfunc{barr}C}=\overline{\limfunc{cone}}\left( \Omega +\limfunc{%
barr}C\right) ,$ Theorem \ref{thm2}\textrm{(i)} reads
\begin{equation*}
\overline{\limfunc{epi}f^{\ast }+K}=\overline{\func{dom}f^{\ast }+\limfunc{%
cone}\left( \Omega \cup \limfunc{barr}C\right) }\times \mathbb{R}
\end{equation*}%
or, equivalently,
\begin{equation*}
\overline{\limfunc{epi}f^{\ast }+K}=\overline{\func{dom}f^{\ast }+\limfunc{%
cone}\left( \Omega +\limfunc{barr}C\right) }\times \mathbb{R}.
\end{equation*}
\end{remark}

The next result provides a sufficient condition to have
\begin{equation*}
\limfunc{epi} (f+\delta_A)^\ast \ \not= \ \overline{\limfunc{epi}f^{\ast }+K}
\end{equation*}
when $A\cap \limfunc{dom} f = \emptyset, A \not= \emptyset$, in terms of $%
A^{\infty}$ and $f^{\infty}$.

\begin{theorem}[On the failure of the hidden assumption with or without
consistency]
\label{thm3}Assume that $A\cap \func{dom}f=\emptyset $ and let $x^{\ast }\in
X^{\ast }.$\ The following statements are equivalent:\newline
$\mathrm{(i)}$ $\limfunc{epi}\left( f+\delta _{A}\right) ^{\ast }\cap \left(
\left\{ x^{\ast }\right\} \times \mathbb{R}\right) \neq \overline{\limfunc{%
epi}f^{\ast }+K}\cap \left( \left\{ x^{\ast }\right\} \times \mathbb{R}%
\right) .$\newline
$\mathrm{(ii)}$ $\left\{ x^{\ast }\right\} \times \mathbb{R}\ \nsubseteq \
\overline{\limfunc{epi}f^{\ast }+K}.$\newline
$\mathrm{(iii)}$ $x^{\ast }\notin \overline{\func{dom}f^{\ast }+\limfunc{cone%
}\Omega +\limfunc{barr}C}.$\newline
$\mathrm{(iv)}$ There exists $d\in \bigcap\limits_{i\in I}\left[
f_{i}^{\infty }\leq 0\right] \cap C^{\infty }$ such that $f^{\infty }\left(
d\right) <\left\langle x^{\ast },d\right\rangle .$\newline
Moreover, if $A\neq \emptyset ,$ then we can add\newline
$\mathrm{(v)}$ There exists $d\in A^{\infty }$ such that $f^{\infty }\left(
d\right) <\left\langle x^{\ast },d\right\rangle .$\newline
\end{theorem}

%\red{+++++++++++}

\textbf{Proof.} \ Since $A\cap \func{dom}f=\emptyset $ we have $\limfunc{epi}%
\left( f+\delta _{A}\right) ^{\ast }=X^{\ast }\times \mathbb{R}$ and the
equivalence $\left[ \mathrm{(i)}\Longleftrightarrow \mathrm{(ii)}\right] $
is obvious.

By Theorem \ref{thm2} we have%
\begin{equation*}
\overline{\limfunc{epi}f^{\ast }+K}=\overline{\func{dom}f^{\ast }+\limfunc{%
cone}\Omega +\limfunc{barr}C}\times \mathbb{R}
\end{equation*}%
and \textrm{(ii)} means
\begin{equation*}
\left\{ x^{\ast }\right\} \times \mathbb{R\nsubseteq }\overline{\func{dom}%
f^{\ast }+\limfunc{cone}\Omega +\limfunc{barr}C}\times \mathbb{R}
\end{equation*}%
or, equivalently,
\begin{equation*}
x^{\ast }\notin \overline{\func{dom}f^{\ast }+\limfunc{cone}\Omega +\limfunc{%
barr}C}.
\end{equation*}%
So, $\left[ \mathrm{(ii)}\Longleftrightarrow \mathrm{(iii)}\right] .$

$\left[ \mathrm{(iii)}\Longleftrightarrow \mathrm{(iv)}\right] $ By Lemma %
\ref{lem1nw}, (\ref{2.3}), (\ref{2.7}) and (\ref{4.2}), there is $d\in
X\diagdown \left\{ 0_{X}\right\} $ such that
\begin{equation}
\begin{array}{ll}
\left\langle x^{\ast },d\right\rangle & >\delta _{\left( \func{dom}f^{\ast }+%
\limfunc{cone}\Omega +\limfunc{barr}C\right) }^{\ast }\left( d\right) \\
& =\delta _{\func{dom}f^{\ast }}^{\ast }\left( d\right) +\delta _{\limfunc{%
cone}\Omega }^{\ast }\left( d\right) +\delta _{\limfunc{barr}C}^{\ast
}\left( d\right) \\
& =f^{\infty }\left( d\right) +\delta _{\Omega ^{-}}\left( d\right) +\delta
_{C^{\infty }}\left( d\right) .%
\end{array}
\label{4.7}
\end{equation}

Since, by (\ref{2.8}),
\begin{equation*}
\Omega ^{-}=\left( \bigcup\limits_{i\in I}\func{dom}f_{\QTR{frametitle}{i}%
}^{\ast }\right) ^{-}=\bigcap\limits_{\QTR{frametitle}{i}\in
\QTR{frametitle}{I}}\left( \func{dom}f_{\QTR{frametitle}{i}}^{\ast }\right)
^{-}=\bigcap\limits_{\QTR{frametitle}{i}\in \QTR{frametitle}{I}}\left[
f_{i}^{\infty }\leq 0\right] ,
\end{equation*}%
we obtain from (\ref{4.7}) that \textrm{(iii)} is equivalent to the
existence of $d\in X$ such that $\left\langle x^{\ast },d\right\rangle
>f^{\infty }\left( d\right) $ and $d\in \bigcap\limits_{i\in I}\left[
f_{i}^{\infty }\leq 0\right] \cap C^{\infty },$ {which is \textrm{(iv)}. }

If, moreover, $A\neq \emptyset ,$\ by (\ref{2.9}) we have%
\begin{equation*}
\begin{array}{ll}
A^{\infty }=B^{\infty }\cap C^{\infty } & =\left( \bigcap\limits_{%
\QTR{frametitle}{i}\in \QTR{frametitle}{I}}\left[ f_{i}\leq 0\right] \right)
^{\infty }\cap C^{\infty } \\
& =\bigcap\limits_{\QTR{frametitle}{i}\in \QTR{frametitle}{I}}\left[
f_{i}\leq 0\right] ^{\infty }\cap C^{\infty } \\
& =\bigcap\limits_{\QTR{frametitle}{i}\in \QTR{frametitle}{I}}\left[
f_{i}^{\infty }\leq 0\right] \cap C^{\infty }%
\end{array}%
\end{equation*}%
and $\left[ \mathrm{(iv)}\Longleftrightarrow \mathrm{(v)}\right] $ holds.$%
\hfill $\noindent $\hfill $\noindent $\square $

\begin{remark}
\label{rem2}In the counterexample (Example \ref{Example1}), condition
\textrm{(v)} in Theorem \ref{thm3} is satisfied with $x^{\ast }=0_{X^{\ast
}} $ \thinspace and $d=\left( 1,1\right) .$
\end{remark}

\section{Applications}

We firstly go a bit further on the consistency/inconsistency of the system $%
\sigma $ in (\ref{1.2}).

\begin{proposition}[The closure of the characteristic cone of inconsistent
systems]
\label{prop1}The following statements are equivalent:\newline
$\mathrm{(i)}$ $\sigma $ is inconsistent.\newline
$\mathrm{(ii)}$ $\overline{K}=\overline{\limfunc{cone}\Omega +\limfunc{barr}C%
}\times \mathbb{R}.$\newline
$\mathrm{(iii)}$ $\overline{K}=\overline{\limfunc{cone}}\left( \Omega \cup
\limfunc{barr}C\right) \times \mathbb{R}.$\newline
$\mathrm{(iv)}$ $\overline{K}=\overline{\limfunc{cone}}\left( \Omega +%
\limfunc{barr}C\right) \times \mathbb{R}.$
\end{proposition}

\textbf{Proof.} \ Apply Theorem \ref{thm2} to the case when $f\equiv 0.$ We
have $f^{\ast }=\delta _{\left\{ 0_{X^{\ast }}\right\} },$ $\limfunc{epi}%
f^{\ast }=\left\{ 0_{X^{\ast }}\right\} \times \mathbb{R}_{+}$ and $\func{dom%
}f^{\ast }=\left\{ 0_{X^{\ast }}\right\} .$ To conclude the proof, note that
the characteristic cone $K$ has the epigraphic property, $K+\left\{
0_{X^{\ast }}\right\} \times \mathbb{R}_{+}=K,$ and take Remark \ref{rem1}
into account.$\hfill $\noindent $\hfill $\noindent $\square \medskip $

We now recover a classical result about the consistency of infinite convex
systems \cite[Theorem 3.1]{DGL06}.

\begin{corollary}[Existence theorem]
\label{corol2}$\text{\textrm{\cite[Theorem\ 3.1]{DGL06}}}$ The following
statements are equivalent to each other:\newline
$\mathrm{(a)}$ $\sigma $ is consistent. \newline
$\mathrm{(b)}$ $\left( 0_{X^{\ast }},-1\right) \notin \mathrm{\,}\overline{%
\limfunc{cone}}\left\{ \bigcup\limits_{i\in I}\func{gph}f_{i}^{\ast }\cup
\func{gph}\delta _{C}^{\ast }\right\} .$\newline
$\mathrm{(c)}$ $\left( 0_{X^{\ast }},-1\right) \notin \overline{\limfunc{cone%
}}\left\{ \limfunc{epi}\delta _{C}^{\ast }\cup \bigcup_{\QTR{frametitle}{i}%
\in I}\limfunc{epi}f_{\QTR{frametitle}{i}}^{\ast }\right\} .$ \newline
$\mathrm{(d)}$ $\overline{\limfunc{cone}}\left\{ \limfunc{epi}\delta
_{C}^{\ast }\cup \bigcup_{\QTR{frametitle}{i}\in I}\limfunc{epi}f_{%
\QTR{frametitle}{i}}^{\ast }\right\} \neq \mathrm{\,}\overline{\limfunc{cone}%
}\left\{ \bigcup_{i\in I}\func{dom}f_{i}^{\ast }\cup \func{dom}\delta
_{C}^{\ast }\right\} \times \mathbb{R}.$
\end{corollary}

\textbf{Proof.} Statements $\rceil \mathrm{(a)}$ and $\rceil \mathrm{(d)}$
coincide with statement $\mathrm{(i)}$ and $\mathrm{(iv)}$ of Proposition %
\ref{prop1}, respectively, while the equivalences $\rceil \mathrm{(b)}%
\Leftrightarrow \rceil \mathrm{(c)}\Leftrightarrow \rceil \mathrm{(d)}$ can
easily be obtained via convex analysis arguments.$\hfill $\noindent $\hfill $%
\noindent $\square \medskip $

We finally associate with the given data, $f\in \Gamma \left( X\right) $, $C$
and $\left\{ f_{i}\right\} _{i\in I}$\ as in (\ref{1.2}), the following
linearly perturbed optimization problem:
\begin{equation*}
\begin{array}{lll}
{(P_{x^{\ast }})}: & \limfunc{Min}_{x\in C} & f\left( x\right) -\left\langle
x^{\ast },x\right\rangle \\
& \text{subject to} & f_{i}(x)\leq 0,i\in I,%
\end{array}%
\text{ }
\end{equation*}%
where $x^{\ast }\in X^{\ast }.$

We now apply the results obtained in the previous section to derive
optimality conditions (both in asymptotic and non-asymptotic forms) for the
linearly perturbed problem ${(P_{x^{\ast }})}$.

\begin{proposition}[An assymptotic optimality condition]
\label{prop2}Let $\overline{x}\in A\cap \func{dom}f$ and $x^{\ast }\in
X^{\ast }.$ Then, $\overline{x}$ is an optimal solution of $(P_{x^{\ast }})$
if and only if
\begin{equation*}
\left( x^{\ast }, \left\langle x^{\ast },\overline{x}\right\rangle -f\left(
\overline{x}\right) \right) \in\ \overline{\limfunc{epi}f^{\ast }+K}.
\end{equation*}
\end{proposition}

\textbf{Proof.} The optimality of $\bar{x}$ for ${(P_{x^{\ast }})}$ means
that
\begin{equation*}
x\in A\Longrightarrow f\left( x\right) -\left\langle x^{\ast
},x\right\rangle \geq f\left( \overline{x}\right) -\left\langle x^{\ast },%
\overline{x}\right\rangle .
\end{equation*}

Applying Theorem \ref{thm1} with $s=f\left( \overline{x}\right)
-\left\langle x^{\ast },\overline{x}\right\rangle $, this holds if and only
if
\begin{equation*}
\Big(x^{\ast },\left\langle x^{\ast },\overline{x}\right\rangle -f\left(
\overline{x}\right) \Big)\in \overline{\limfunc{epi}f^{\ast }+K}.\ \ \ \ \ \
\ \ \ \ \ \ \ \ \hfill \square
\end{equation*}

\begin{proposition}[KKT characterization]
\label{prop4} Assume that $\limfunc{epi}f^{\ast }+K$ is $w^{\ast }-$closed
and let $x^{\ast }\in X^{\ast }$ and $\overline{x}\in A\cap \func{dom}f$.
The following statements are equivalent:\newline
$\mathrm{(i)}$ $\overline{x}$ is an optimal solution of $(P_{x^{\ast }}).$%
\newline
$\mathrm{(ii)}$ There exist $J\in \mathcal{J}$ and $\lambda _{j}\geq 0$ for
all $j\in J$ such that
\begin{eqnarray*}
x^{\ast }\!\!\! &\in &\!\!\!\partial f(\overline{x})+N_{C}(\overline{x}%
)+\sum_{j\in J}\lambda _{j}\partial f_{j}(\overline{x}), \\
&&\lambda _{j}f_{j}(\overline{x})=0,\forall j\in J.
\end{eqnarray*}%
\newline
$\mathrm{(iii)}$ There exist $J\in \mathcal{J}$ and $\lambda _{j}\geq 0$ for
all $j\in J$ such that%
\begin{eqnarray*}
x^{\ast }\!\!\! &\in &\!\!\!\partial \Big( f+\delta _{C}+\sum_{j\in
J}\lambda _{j} f_{j}\Big)(\overline{x}), \\
&&\lambda _{j}f_{j}(\overline{x})=0,\forall j\in J.
\end{eqnarray*}
\end{proposition}

\textbf{Proof.} $[\mathrm{(i)}\Longrightarrow \mathrm{(ii)}]$ Let $\overline{%
x}\in A\cap \func{dom}f$ be an optimal solution of $(P_{x^{\ast }})$.
Applying Theorem \ref{thm1b} with $s=f\left( \overline{x}\right)
-\left\langle x^{\ast },\overline{x}\right\rangle $, there exist $J\in
\mathcal{J}$, $\lambda _{j}\geq 0$, $u_{j}^{\ast }\in \limfunc{dom}%
f_{j}^{\ast }$ for all $j\in J$, $u^{\ast }\in \limfunc{dom}f^{\ast }$ and $%
v^{\ast }\in \limfunc{barr}C$ such that $x^{\ast }=u^{\ast }+v^{\ast
}\sum_{j\in J}\lambda _{j}u_{j}^{\ast }$ and
\begin{equation*}
\begin{array}{ll}
0\geq & \left[ f^{\ast }\left( u^{\ast }\right) +f\left( \overline{x}\right)
-\left\langle u^{\ast },\overline{x}\right\rangle \right] +\left[ \delta
_{C}^{\ast }\left( v^{\ast }\right) -\left\langle v^{\ast },\overline{x}%
\right\rangle \right] \\
& +\sum_{j\in J}\left[ \lambda _{j}\left( f_{j}^{\ast }\left( u_{j}^{\ast
}\right) +f_{j}\left( \overline{x}\right) -\left\langle u_{j}^{\ast },%
\overline{x}\right\rangle \right) \right] +\sum_{j\in J}\left[ -\lambda
_{j}f_{j}\left( \overline{x}\right) \right] .%
\end{array}%
\end{equation*}%
Since all the brackets are non-negative, they are all equal to zero, that
entails $u^{\ast }\in \partial f\left( \overline{x}\right) ,$ $v^{\ast }\in
N_{C}(\overline{x}),$ $u_{j}^{\ast }\in \partial f_{j}\left( \overline{x}%
\right) ,$ for all $j\in J_{0}:=\left\{ j\in J:\lambda _{j}>0\right\} ,$ and
$\lambda _{j}f_{j}(\overline{x})=0,$ for all $j\in J.$ We then have under
the rule that the empty sum is $0_{X^{\ast }},$%
\begin{equation*}
\begin{array}{ll}
x^{\ast } & \in \ \ \partial f(\overline{x})+N_{C}(\overline{x})+\sum_{j\in
J_{0}}\lambda _{j}\partial f_{j}(\overline{x}), \\
& =\ \ \partial f(\overline{x})+N_{C}(\overline{x})+\sum_{j\in J}\lambda
_{j}\partial f_{j}(\overline{x})%
\end{array}%
\end{equation*}%
and $\mathrm{(ii)}$ has been proved.

$[\mathrm{(ii)}\Longrightarrow \mathrm{(iii)}]$ This is obvious as
\begin{equation*}
\partial f(\overline{x})+\sum_{j\in J}\lambda _{j}\partial f_{j}(\overline{x}%
)+N_{C}(\overline{x})\subset \partial \Big(f+\sum_{j\in J}\lambda
_{j}f_{j}+\delta _{C}\Big)(\overline{x}).
\end{equation*}

$[\mathrm{(iii)}\Longrightarrow \mathrm{(i)}]$ If $\mathrm{(iii)}$ holds
then
\begin{equation*}
\Big(f+\sum_{j\in J}\lambda _{j}f_{j}+\delta _{C}\Big)(x)-\Big(f+\sum_{j\in
J}\lambda _{j}f_{j}+\delta _{C}\Big)(\overline{x})\geq \langle x^{\ast },x-%
\overline{x}\rangle ,
\end{equation*}%
or, equivalently (as $\sum_{j\in J}\lambda _{j}f_{j}(\overline{x})=0$),
\begin{equation*}
f(x)-\langle x^{\ast },x\rangle +\sum_{j\in J}\lambda _{j}f_{j}(x)\geq f(%
\overline{x})-\langle x^{\ast },\overline{x}\rangle +\sum_{j\in J}\lambda
_{j}f_{j}(\overline{x})=f(\overline{x})-\langle x^{\ast },\overline{x}%
\rangle ,\forall x\in C.
\end{equation*}%
Thus, when $x\in A$ we get
\begin{equation*}
f(x)-\langle x^{\ast },x\rangle \geq f(\overline{x})-\langle x^{\ast },%
\overline{x}\rangle ,
\end{equation*}%
which means that $\overline{x}$ is an optimal solution of $(P_{x^{\ast }})$
and $\mathrm{(i)}$ has been proved. The proof is complete. \hfill $\square $

\bigskip

\textbf{Disclosure statement}

No potential conflict of interest was reported by the author(s).$\bigskip $

\textbf{Funding}$\bigskip $

\textbf{Acknowledgements}

The authors are grateful to Marco A. L\'{o}pez for his valuable comments and
suggestions on the manuscript, as well as for his private communication that
a different proof of Corollary \ref{corol1} (assuming also that $A\cap \func{%
dom}f\neq \emptyset $) will appear independently in his forthcoming book
with Rafael Correa and Abderrahim Hantoute \cite{CHL23}.$\bigskip $

\end{document}